\def\Bbb#1{{\bf #1}}

\def\fnote#1{\footnote}
\def\blacksquare{\hbox{\vrule width 4pt height 4pt depth 0pt}}

\def\cwleftpar#1#2{\leftskip #1 \rightskip #2 plus 1fill}
\def\cwrightpar#1#2{\leftskip #1 plus 1fill \rightskip #2}
\def\cwcenterpar#1#2{\leftskip #1 plus 1fill \rightskip #2 plus 1fill}
\def\cwfullpar#1#2{\leftskip#1\rightskip#2}

\def\cwoutdent#1#2{\llap{\hbox to #1{#2 \hss}}\ignorespaces}
\def\cwparbegin#1#2#3#4#5{
	\ifcase #1 \cwleftpar{#2}{#3}
	\or \cwrightpar{#2}{#3}
	\or \cwcenterpar{#2}{#3}
	\else \cwfullpar{#2}{#3}\fi
	\ifcase #4 \baselineskip = 1.5\baselineskip
	\or \baselineskip = 2\baselineskip
	\or \baselineskip = 3\baselineskip
	\else \baselineskip = 1\baselineskip\fi
	\ifdim #5 > 0in \else \noindent \fi
	\noindent\ignorespaces}
\documentclass{article}
\begin{document}
%------------------------------------------------------------------------
%                               ChiWriter FOOTER
%------------------------------------------------------------------------
\advance \vsize by -1\baselineskip
\def\makefootline{
{\vskip \baselineskip \noindent \folio                               \par
}}
%------------------------------------------------------------------------
%                               Your Document
%------------------------------------------------------------------------

\vspace*{2ex}
\noindent {\Huge Linear Transports along Paths in\\[0.4ex] Vector
	Bundles}\\[1.8ex]
\noindent {\Large II. Some Applications}
\vspace*{2ex}

\medskip \noindent Bozhidar Zakhariev Iliev
\fnote{0}{\noindent $^{\hbox{}}$Permanent address:
Laboratory of Mathematical Modeling in Physics,
Institute for Nuclear Research and \mbox{Nuclear} Energy,
Bulgarian Academy of Sciences,
Boul.\ Tzarigradsko chauss\'ee~72, 1784 Sofia, Bulgaria\\
\indent E-mail address: bozho@inrne.bas.bg\\
\indent URL: http://theo.inrne.bas.bg/$\sim$bozho/}

\vspace*{2ex}
{\bf \noindent Published: Communication JINR, E5-93-260, Dubna, 1993}\\[1ex]
\hphantom{\bf Published: }
http://www.arXiv.org e-Print archive No.~math.DG/0412010\\[2ex]

2003 PACS numbers: 02.40.Ma, 02.40.Vh, 04.90.+e\\[2ex]

\noindent
{\small
The \LaTeXe\ source file of this paper was produced by converting a
ChiWriter 3.16 source file into
ChiWriter 4.0 file and then converting the latter file into a
\LaTeX\ 2.09 source file, which was manually edited for correcting numerous
errors and for improving the appearance of the text.  As a result of this
procedure, some errors in the text may exist.
}\\[2ex]

	\begin{abstract}
The linear transports along paths in vector bundles introduced in Ref. [1]
are applied to the special case of tensor bundles over a given differentiable
manifold. Links with the transports along paths generated by derivations
of tensor algebras are investigated. A possible generalization of the theory
of geodesics is proposed when the parallel transport generated by a linear
connection is replaced with an arbitrary linear transport along paths in the
tangent bundle.
	\end{abstract}\vspace{3ex}

{\bf INTRODUCTION}

\medskip
This work is devoted to some simple applications of the considered in [1] linear transports along paths in vector bundles. It is organized as follows. Sect. 2 contains a collection of the main definitions and results of [1] which will be used here. In Sect. 3 the general theory of linear transports along paths is specialized in the case of tensor bundles over a given differentiable manifold. In particular, the connections with the transports along paths generated by derivations of tensor algebras are investigated. The problem for holonomicity of special bases in which the matrix of a linear transport along paths in a tensor bundle is unit is considered in Sect. 4. In Sect. 5 the concept of linearly transported along paths sections of a vector bundle is introduced. Some properties of these section, as well as their connections with derivations along paths, are investigated. In Sect. 6 the firs steps of one possible generalization of the theory of geodesics paths (curves) is proposed. Here we consider paths in a manifold the tangent vector to which after a linear transport along themselves remains such. Sect. 7 contains remarks concerning special bases for parallel transports generated by linear connections.

Part of the results of this work were found in [2] in the special case of the tangent bundle to a given differentiable manifold.

\medskip
\medskip
 {\bf 2. SUMMARY OF SOME RESULTS}

\medskip
In this section we summarize certain needed for this paper definitions and results of [1].

A  linear transport $(L$-transport$) L^{\gamma }_{s  \to t}$along the path
$\gamma :J  \to B, J$ being a real interval, from $s$ to $t, s,t\in J$ in the
vector bundle $(E,\pi ,B)$ with a base $B$, total bundle space $E$ and
projection $\pi :E  \to B$, has the properties:
\[
 L^{\gamma }_{s  \to t}:\pi ^{-1}(\gamma (s))  \to \pi ^{-1}(\gamma (t)),
\qquad (2.1)
\]
\[
 L^{\gamma }_{s  \to t}(\lambda u+\mu v)=\lambda S^{\gamma
}_{s  \to t}u+\mu S^{\gamma }_{s  \to t}v,
\quad \lambda ,\mu \in {\Bbb R},\ u,v\in \pi ^{-1}(\gamma (s))  \qquad (2.2)
\]
\[
  L^{\gamma }_{t  \to r}\circ L^{\gamma }_{s  \to
t}=L^{\gamma }_{s  \to r},\quad  r,s,t\in J,\qquad (2.3)
\]
\[
 L^{\gamma }_{s  \to
s}={\it id}_{\pi ^{-1_{(\gamma (s))}}},\qquad (2.4)
\]
 where {\it id}$_{X}$ means the identity map of the set X. Its general form is
\[
 L^{\gamma }_{s  \to t}={\bigl(}F^{\gamma }_{t}{\bigr)}^{-1}\circ
{\bigl(}F^{\gamma }_{s}{\bigr)},\quad  s,t\in J. \qquad (2.5)
\]
 where $F^{\gamma }_{s}:\pi ^{-1}(\gamma (s))  \to V, V$ being $a \dim(\pi
^{-1}(x))$-dimensional, $x\in B$, vector space, are linear isomorphisms.

 If $\{e_{i}, i=1,\ldots  , \dim(\pi ^{-1}(x)), x\in B\}$ is a field of bases
(or simply a basis) along $\gamma $, i.e. $\{e_{i}(s)\}$ is a basis in $\pi
^{-1}(\gamma (s))$ for $s\in J$, then the matrix $H:(t,s;\gamma )\mapsto
H(t,s;\gamma ):=H^{i}_{.j}(t,s;\gamma )$, $s,t\in J$ of the transport is
defined by
\[
 L^{\gamma }_{s  \to t}e_{i}(s)=H^{j}_{.i}(t,s;\gamma )e_{j}(t),
\quad s,t\in J,\qquad (2.6)
\]
 where here and from here on in our text the Latin indices run
from 1 to $\dim(\pi ^{-1}(x)), x\in B$ and the usual summation rule from 1 to
$\dim(\pi ^{-1}(x))$ over repeated on different levels indices is assumed.

  The matrix function $H$ satisfies the equations
\[
 H(r,t;\gamma )H(t,s;\gamma )=H(r,s;\gamma ),\quad  r,s,t\in J,\qquad (2.7)
\]
\[
  H(s,s;\gamma )={\Bbb I}:=diag(1,\ldots  ,1)=:  \delta ^{i}_{j}  ,
\quad s\in
J\qquad (2.8)
\]
and its general form is
\[
 H(t,s;\gamma )=(F(t;\gamma ))^{-1}F(s;\gamma ), s,t\in J,\qquad (2.9)
\]
 where $F(s;\gamma )$ is a nondegenerate matrix function.

 Let $L$ be a smooth, of class $C^{1}$, linear transport along paths and
$\sigma\in Sec^{1}(\xi|_{\gamma (J)})$, where $\xi \mid _{\gamma
(J)}$is the restriction of the vector bundle $\xi =(E,\pi ,B)$ on the set
$\gamma (J)$ and by Sec$^{k}(\xi ) ($resp. Sec$(\xi ))$ is denoted the set of
$C^{k}($resp. all) section over $\xi $. The generated by $L$ derivation along
$\gamma :J  \to B$ is a map ${\cal D}^{\gamma }$:Sec$^{1}(\xi \mid _{\gamma
(J)})  \to $Sec$(\xi \mid _{\gamma (J)})$ defined by $(s,s+\epsilon \in J)$
\[
{\bigl(}{\cal D}^{\gamma }\sigma {\bigr)}(\gamma (s))
:={\cal D}^{\gamma }_{s}\sigma
:=\lim_{\varepsilon\to 0}
\bigl[ \frac{1}{\varepsilon} {\bigl(}
L^{\gamma }_{s+\epsilon \to s}\sigma (\gamma (s+\epsilon ))
- \sigma (\gamma (s)){\bigr)}\bigr] \qquad(2.10)
\]
 and has the properties:
\[
 {\cal D}^{\gamma }(\lambda \sigma +\mu \tau )=\lambda {\cal D}^{\gamma
}\sigma +\mu {\cal D}^{\gamma }
\qquad \tau , \lambda ,\mu \in {\Bbb R},\
\sigma ,\tau\in Sec^{1}(\xi \mid _{\gamma (J)}),\qquad (2.11)
\]
\[
 {\cal D}^{\gamma }_{s}(f\cdot \sigma )
=\frac{df(s)}{ds} \cdot \sigma (\gamma (s))
+f(s)\cdot ({\cal D}^{\gamma }_{s}\sigma ),\qquad (2.12)
\]
\[
{\cal D}^{\gamma }_{t}\circ L^{\gamma }_{s  \to t}\equiv 0, \quad s,t\in
J\qquad (2.13)
\]
 for any $C^{1}$function $f:J  \to {\Bbb R}$.

 The coefficients $\Gamma ^{i}_{.j}(s;\gamma )$ of $L$ along $\gamma $ in
$\{e_{i}\}$ coincide with those of ${\cal D}^{\gamma }$and are given, e.g.,
by
\[
{\cal D}^{\gamma }_{s}e_{j}=\Gamma ^{i}_{.j}(s;\gamma )e_{i}(s).\qquad
(2.14)
\]
 The explicit connection of $\Gamma ^{i}_{.j}(s;\gamma )$ with the
matrix $H$ of $L$ is
\[
 \Gamma _{\gamma }(s)
= [ \Gamma ^{i}_{.j}(s;\gamma ) ]
= \frac{\partial H(s,t;\gamma)}{\partial t} \Big|_{t=s}
= F^{-1}(s;\gamma )  \frac{d F(s;\gamma)}{ds}. \qquad (2.15)
\]
 If $\sigma =\sigma ^{i}e_{i}\in $Sec$^{1}(\xi \mid _{\gamma (J)})$, then explicitly
(2.10) reads
\[
{\bigl(}{\cal D}^{\gamma }\sigma {\bigr)}(\gamma (s))
= \Bigl[ \frac{d \sigma^i(\gamma(s))}{ds}
+ \Gamma ^{i}_{.j}(s;\gamma )\sigma ^{j}(\gamma (s))
\Bigr]  e_{i}(s). \qquad (2.16)
\]

\newpage

 {\bf 3. LINEAR TRANSPORTS ALONG PATHS IN\\ TENSOR BUNDLES}

\medskip
Now we shall consider linear transports along paths in the special case of
tensor bundles over a given differentiable manifold. In particular we shall
investigate the ties of these transports with the ones generated by
derivations of tensor algebras (the $S$-transports).

Let $M$ be a differentiable manifold,
$T^{p}_{.q}\mid _{x}(M),$ $p,q\in {\Bbb N}\cup \{0\}$ be the tensor space of
type   q   over $M$ at $x\in M$, and
$(T^{p}_{.q}(M),\pi ,M)$,
$T^{p}_{.q}(M):=\cup_{x\in M} T^{p}_{.q}|_{x}$ be
the tensor bundle of type q   over $M$ with a projection $\pi :T^{p}_{.q}(M)
\to M$ such that $\pi ^{-1}(x):=T^{p}_{.q}\mid _{x}(M), x\in M (cf. [3])$.

Till the end of this section we shall deal with $L$-transports
$^{p}_{q}L^{\gamma }$along paths $\gamma :J\to M$ acting, respectively, in
the tensor bundles $(T^{p}_{.q}(M),\pi ,M)$ of arbitrary type   q  , i.e. we
will investigate maps
\[
 ^{p}_{q}L:\gamma \mapsto ^{p}_{q}L^{\gamma }:(s,t)\mapsto
^{p}_{q}L^{\gamma }_{s\to t}:T^{p}_{.q}|_{\gamma (s)}(M)\to
		T^{p}_{.q}|_{\gamma (t)}(M)  \qquad (3.1)
\]
 satisfying $(2.2)-(2.4)$.

${\bf P}$ractically is far  more convenient instead of the transports
$^{p}_{q}L^{\gamma }, p,q\ge 0$ along $\gamma $ to be used the equivalent to
them map $L^{\gamma }$, {\it the} $(L-)${\it transport along} $\gamma $ {\it
in the tensor algebra over} $M$, defined by $L^{\gamma }:(s,t)\mapsto
L^{\gamma }_{s\to t}, s,t\in J$, where
\[
L^{\gamma }_{s\to t}: \bigcup_{p,q=0}^{\infty} T^{p}_{.q}|_{\gamma (s)}(M)
 \to \bigcup_{p,q=0}^{\infty} T^{p}_{.q}|_{\gamma (t)}(M),
\quad  s,t\in J \qquad (3.2)
\]
 with
\[
L^{\gamma }_{s\to t}T_{0}
:=^{p}_{q}L^{\gamma }_{s\to t}T_{0} \quad
for\ T_{0}\in T^{p}_{.q}|_{\gamma (s)}(M).\qquad (3.3)
\]
 Before formulating certain results concerning the so
defined $L$-transports $L^{\gamma }$along $\gamma $ we will present the
concrete form of some formulae from the preceding section in the case of the
maps (3.2).

According to the equality (2.6) the matrix elements
$H^{...}_{...}(t,s;\gamma )$ of $L^{\gamma }$are uniquely defined by the
expansion
\[
L^{\gamma }_{s\to t}{\bigl(}
E^{j_{1}\ldots j_q}_{i_{1}\ldots i_p}|_{\gamma (s)}{\bigr)}
=H^{k_{1}\ldots k_p;j_1\ldots j_q}
  _{l_{1}\ldots l_q;i_1\ldots i_p}(t,s;\gamma )
E^{l_{1}\ldots l_q}_{k_{1}\ldots k_p}|_{\gamma (t)}, \qquad (3.4)
\]
 and are components of a two-point tensor from $T^{p}_{.q}$  $_{\gamma (t)}(M)\otimes
\otimes T^{q}_{.p}$  $_{\gamma (s)}(M)$. Here
 $E^{l_{1}\ldots l_q}_{k_{1}\ldots k_p}
:=E^{l_{1}}\otimes \cdot \cdot \cdot E^{l_{q}}\otimes
E_{k_{1}}\otimes \cdot \cdot \cdot \otimes E_{k_{p}}, \otimes $ being the
tensor product sign, is a basis (field of bases) in $T^{p}_{.q}(M)$ generated
by the bases $\{E^{i}\}$ and its dual $\{E_{j}\}$ in, respectively,
$T^{1}_{.0}(M)$ and $T^{0}_{.1}(M)$.

If
\(
T=T^{j_{1}\ldots j_q}_{i_{1}\ldots i_p}
	E^{i_{1}\ldots i_p}_{j_{1}\ldots j_q}
|_{\gamma (s)}\in T^{p}_{.q}|_{\gamma (s)}(M),
\)
then, due to (2.6) and the linearity of $L^{\gamma }_{s\to t}$, we have
\[
L^{\gamma }_{s\to t}T
=
H^{k_{1}\ldots k_p;j_{1}\ldots j_q}
 _{l_{1}\ldots l_q;i_{1}\ldots i_p}(t,s;\gamma )
	T^{i_{1}\ldots i_p}_{j_{1}\ldots j_q}
\Bigl( E^{l_{1}\ldots l_q}_{k_{1}\ldots k_p} |_{\gamma (t)}\Bigr).
 \qquad (3.5)
\]
Because of (2.7) and (2.8) the following equalities are valid:
\[
H^{k_{1}\ldots k_p;j_{1}\ldots j_q}
 _{l_{1}\ldots l_q;i_{1}\ldots i_p}(r,t;\gamma)
\cdot
H^{i_{1}\ldots i_p;n_{1}\ldots n_q}
 _{j_{1}\ldots j_q;m_{1}\ldots m_p}(t,s;\gamma )=
\]
\[
=
H^{k_{1}\ldots k_p;n_{1}\ldots n_q}
 _{l_{1}\ldots l_q;m_{1}\ldots m_p}(r,s;\gamma ),\quad  r,s,t\in J,
\qquad (3.6)
\]
\[
 H^{k_{1}\ldots k_p;j_{1}\ldots j_q}
  _{l_{1}\ldots l_q;i_{1}\ldots i_p}(s,s;\gamma )
=\Bigl( \prod_{a=1}^{p} \delta ^{k_{a}}_{i_{a}} \Bigr)
 \Bigl( \prod_{b=1}^{q} \delta ^{j_{b}}_{l_{b}} \Bigr),\qquad (3.7)
\]
 which, in our case, are equivalent to (2.3) and (2.4) respectively.

From (3.6) and (3.7) it can easily be obtained the general form of the matrix elements $H^{...}_{...}(\ldots  )$. Its explicit form is described by the component form of (2.9) in which the indexes must be replaced with the corresponding multi-indexes (e.g. $i\mapsto (i_{1},\ldots  ,i_{p})$ etc.).

{\bf Definition 3.1.} The $L$-transport $L$ in the tensor algebra over $M$
will be called consistent (resp. along $\gamma )$ with the operation tensor
multiplication if
\[
 L^{\gamma }_{s\to t}(A\otimes B)=(L^{\gamma }_{s\to t}A)\otimes
(L^{\gamma }_{s\to t}B), s,t\in J\qquad (3.8)
\]
 holds for arbitrary tensors A
and $B$ at the point $\gamma (s)$ and any (resp. the given) path $\gamma $.

 It is easily seen (3.8) to be equivalent to
\[
  ^{p_{1}}_{q_{1}}L^{\gamma }_{s\to t}={\bigl(}^{p_{1}}_{q_{1}}L^{\gamma
}_{s\to t}{\bigr)}\otimes {\bigl(}^{p_{2}}_{q_{2}}L^{\gamma }_{s\to
t}{\bigr)}, s,t\in J\qquad (3.9)
\]
 for arbitrary nonnegative integers $p_{1}, p_{2}, q_{1}$ and $q_{2}$.

With this we end the necessary for the following preliminary material.

{\bf Proposition 3.1.} The equality (3.8) is fulfilled iff the defined by
(3.4) matrix elements of $L^{\gamma }$have the representation
\[
H^{k_{1}\ldots k_p;j_{1}\ldots j_q}
  _{l_{1}\ldots l_q;i_{1}\ldots i_p}(t,s;\gamma )
= \Bigl( \prod_{a=1}^{p} H^{k_{a}}_{..i_{a}}(t,s;\gamma ) \Bigr)
  \Bigl( \prod_{b=1}^{q} H^{..^{j_{b}}}_{l_{b}}(t,s;\gamma ) \Bigr)
 \qquad (3.10)
\]
 where $H^{\hbox{i.}}_{.j}(t,s;\gamma )$ and $H^{.j}_{\hbox{i.}}(t,s;\gamma
)$ define, according to (3.4), the transports along $\gamma $ of the vectors
and covectors respectively.

 {\bf Proof.} (3.10) follows from $(3.4), (p+q-1)$ times application of (3.8)
to the tensor product $A_{1}\otimes \cdot \cdot \cdot \otimes A_{p}\otimes
B_{1}\otimes \cdot \cdot \cdot \otimes B_{q}$for arbitrary $A_{1},\ldots
,A_{p}\in T_{\gamma (s)}(M)$ and $B_{1},\ldots  ,B_{q}\in T^{*}_{\gamma
(s)}(M)$ and the arbitrariness of these vectors and covectors. On the
opposite, if (3.10) is valid, then an elementary check with the help of (3.4)
shows that (3.8) is true.\blacksquare

{\bf Corollary 3.1.} If the $L$-transport $L^{\gamma }$along $\gamma $ in the tensor algebra over $M$ is consistent with the tensor product, then it is uniquely defined if its action is given on vectors and covectors.

{\bf Proof.} This result is a direct consequence from proposition 3.1 and equality (3.5).\blacksquare

{\bf Corollary 3.2.} If the $L$-transport $L^{\gamma }$along $\gamma $ in
the tensor algebra over $M$ is consistent with the tensor product, then
\[
L^{\gamma }_{s\to t}(\lambda )=\lambda ,\quad s,t\in J,\
 \lambda \in {\Bbb R}.\qquad (3.11)
\]

{\bf Proof.} Putting $A=B=1\in {\Bbb R}$ in (3.8), we find
$L^{\gamma }_{s\to t}(1)=1$ which, due to (2.2), is equivalent to (3.11).
(The same result follows also from (3.5) and (3.10) for $p=q=0\in {\Bbb
R}.)\blacksquare $

Let us note that in the general case, due to (2.2), instead of (3.11), we
have
\[
L^{\gamma }_{s\to t}(\lambda )=h(t,s;\gamma )\cdot \lambda ,
\quad  s,t\in J, \ \lambda \in {\Bbb R},\qquad (3.12)
\]
 where (see (2.2)) the 2-point scalar $h(t,s;\gamma )$ is defined by
\[
 h(t,s;\gamma ):=L^{\gamma }_{s\to t}(1),\quad s,t\in J.  \qquad (3.13)
\]
According to (3.6) and (3..7) it has the properties
\[
 h(r,t;\gamma )h(t,s;\gamma )=h(r,s;\gamma ),\quad r,s,t\in J,\qquad (3.14)
\]
\[
h(s,s;\gamma )=1,\quad s\in J\qquad (3.15)
\]
 and because of (2.9) its general form is
\[
 h(t,s;\gamma )=f(s;\gamma )/f(t;\gamma ),\quad s,t\in J,\qquad (3.16)
\]
 where $f$ is nonvanishing scalar function of $s$ and $\gamma $.

Evidently, the transport $^{0}_{0}L^{\gamma }$is consistent with the tensor product iff it is fulfilled (3.11), which is equivalent to $h(t,s;\gamma )=1, s,t\in J$, or all the same $f(s;\gamma )=f_{o}(\gamma )$, i.e. when $f(s;\gamma )$ does not depend on $s\in $J.

An essential role play the $L$-transports along paths in the tensor algebra
over $M$ which commute with the contraction operator $C$, i.e. transports
$L^{\gamma }$along $\gamma $ for which
\[
L^{\gamma }_{s\to t}\circ C-C\circ L^{\gamma }_{s\to t}=0, s,t\in
J,\qquad (3.17)
\]
 or, written in another way,
\[
{\bigl(}^{p-1}_{q-1}L^{\gamma }_{s\to t}{\bigr)}\circ C-C\circ
{\bigl(}^{p}_{q}L^{\gamma }_{s\to t}{\bigr)}=0, \quad p,q\ge 1,\ s,t\in J.
 \qquad (3.17^\prime )
\]

 {\bf Proposition} ${\bf 3}{\bf .}{\bf 2}{\bf .} A$ given
$L$-transport along $\gamma $ in the tensor algebra over $M$ satisfies
simultaneously (3.8) and (3.17) if and only if its matrix elements are given
by (3.10) in which $H^{\hbox{i.}}_{.j}(t,s;\gamma )$ and
$H^{.j}_{\hbox{i.}}(t,s;\gamma )$, defining the $L$-transports along $\gamma
$ of the vectors and covectors respectively, are elements of mutually inverse
matrices, i.e. they are connected by the relationship
\[
 H^{\hbox{k.}}_{.i}(t,s;\gamma )H^{.j}_{\hbox{k.}}(t,s;\gamma )=\delta
^{\hbox{i.}}_{.j}.\qquad (3.18)
\]

 {\bf Proof.} If (3.10) and (3.18) are true,
then with the help of (3.5), by a direct calculation we confirm ourselves
that (3.8) and (3.17) are identically satisfied. On the opposite, if (3.8)
and (3.17) are true, then by proposition 3.1 is valid (3.10), and applying
(3.17) to the tensor $(E^{j}\mid _{\gamma (s)})\otimes (E_{i}\mid _{\gamma
(s)})$ and using (3.10) for $p=q=1$, we get (3.18).\blacksquare

{\bf Corollary 3.3.} An $L$-transport along paths in the tensor algebra over $M$ which is consistent with the tensor multiplication and commutes with the contraction operator is uniquely defined by fixing its action on vectors or, equivalently, on covectors.

{\bf Proof.} This result follows from (3.5) and propositions 3.1 and 3.2, i.e. from $(3.5), (3.10)$ and (3.18).\blacksquare

{\bf Proposition 3.3.} If a given $L$-transport $L^{\gamma }$along $\gamma
:J\to M$ in the tensor algebra over $M$ satisfies (3.8) and (3.17) and
$\nabla $ is the covariant derivative define by an affine connection with
coefficients $\Gamma ^{i}_{.jk}(x)$ at $x\in M$, then
\[
{\cal D}^{\gamma }_{s}=D_{V}|_{\gamma (s)}, s\in J,\qquad (3.19)
\]
where ${\cal D}^{\gamma }_{s}$ is generated by $L^{\gamma }$ through $(2.10),
V$ is defined on a neighborhood of $\gamma (J)$ vector field with the
property $V_{\gamma (s)}=\dot\gamma(s)$, $s\in J$ and
$D_{V}|_{\gamma (s)}$is the derivative along $\gamma $ at $\gamma (s)$
defined by
\[
D_{V}|_{\gamma (s)}=\nabla _{V}|_{\gamma (s)}+H_{V}|_{\gamma (s)},\qquad (3.20)
\]
 where
\[
 H_{V}|_{\gamma (s)}:=[\nabla _{V}|_{\gamma
(s)}(H^{\hbox{i.}}_{.j}(t,s;\gamma )E_{i}\mid _{\gamma (t)}\otimes E^{j}\mid
_{\gamma (s)})]\mid _{t=s}. \qquad (3.21)
\]
 {\bf Remark.} The restriction $\nabla_{V}\mid _{\gamma (s)}$ means that
$\nabla_{V}$ acts only on the defined at the point $\gamma (s)$ objects, and
consequently (see [3] and (3.7))
\[
\bigl( H_{V}\mid _{\gamma (s)} \bigr)^{\hbox{i.}}_{.j}
= V^{k}_{\gamma (s)}
\Bigl( \frac{\partial H_{.j}^{i.}(t,s;\gamma)}{\partial \gamma^k(s)}
-
\Gamma^{l}_{jk}(\gamma (s))H^{\hbox{i.}}_{.l}(t,s;\gamma )
\Bigr) \Big|_{_{t=s}}=
\]
\[
\qquad
= \frac{\partial H_{.j}^{i.}(t,s;\gamma)}{\partial s} \Big|_{t=s}
- \Gamma ^{i}_{.jk}(\gamma (s)) \dot\gamma^{k}(s).\qquad (3.21^\prime )
\]

{\bf Proof.} If $T$ is a $C^{1}$ tensor field of type $(p,q)$ defined on a
neighborhood of $\gamma(J)$, then due to (2.15) and (2.16), we have
\[
 ({\cal D}^{\gamma }_{s}T)^{i_{1}\ldots i_p}_{j_{1}\ldots j_q}
=\frac{d}{ds}T^{i_{1}\ldots i_p}_{j_{1}\ldots j_q}(\gamma (s)) +
\]
\[
 + \Bigl(
\frac{\partial }{\partial s }
H^{i_{1}\ldots i_p;k_{1}\ldots k_q}
 _{j_{1}\ldots j_q;l_{1}\ldots l_p}(t,s;\gamma )
\Bigr)\Big|_{t=s}
T^{l_{1}\ldots l_p}_{k_{1}\ldots k_q}(\gamma (s)).  \qquad (3.22)
\]

  Substituting here (3.10) (see proposition 3.1), with the help of (3.18)
(see proposition $3.2), (3.21^\prime )$ and the equality (see e.g. [3])
\[
\frac{d}{ds}
T^{i_{1}\ldots i_p}_{j_{1}\ldots j_q}(\gamma (s))
= [(\nabla _{V}T)|_{\gamma (s)}]^{i_{1}\ldots i_p}_{j_{1}\ldots j_q} -
\]
\[
\sum_{a=1}^{p} \Gamma ^{i_{a}}_{..kl}(\gamma (s))
T^{i_{1}\ldots i_{a-1}ki_{a+1}\ldots i_p}_{j_{1}\ldots j_q}(\gamma (s))
\dot\gamma^{l}(s) +
\]
\[
+
\sum_{b=1}^{q} \Gamma ^{k}_{.j_{b}}(\gamma (s))
T^{i_{1}\ldots i_p}_{j_{1}\ldots j_{b-1}ki_{b+1}\ldots j_q}(\gamma (s))
\dot\gamma^{l}(s),
\]
 after an easy transformations, we get ${\cal D}^{\gamma }_{s}T=(\nabla
_{V}+H_{V})\mid _{\gamma (s)}(T(\gamma (s)))$. From here, due to the
arbitrariness of $T$ follows (3.19).\blacksquare

It is clear that taken by itself a consistent with the tensor product and
commuting with the contractions $L$-transport along paths in the tensor
algebra over $M$ does not define globally some derivation, but on a given
path $\gamma $ it uniquely defines the derivation (3.20), the action of
which, due to (3.19), does not depend on the used in its definition
additional covariant derivative $\nabla $. This fact allows us if it is
given an $L$-transport along paths with the above properties and the manifold
$M$ is covered with a $C^{1}$ congruence of paths $\{\gamma _{\lambda }:
\gamma _{\lambda }:J_{\lambda }\to M$,   $_{\lambda   }\gamma _{\lambda
}(J_{\lambda })\left.\begin{array}{c} \end{array}\right) =M, \lambda \in
\Lambda \subset {\Bbb R}^{\dim(M)-1}\}$ to construct a global $S$-derivation
$D_{V}$of the tensor algebra over $M$ in the following way. If $x\in M$, then
there exist unique $\mu (x)\in \Lambda $ and $s^{o}_{\mu (x)}\in J_{\mu
(x)}$such that $x=\gamma _{\mu (x)}(s^{o}_{\mu (x)})$. We define the vector
field $V$ by the equality
$V_{x}:=\dot\gamma_{\mu (x)}(s^{o}_{\mu (x)}), x\in M$ and put
\[
 D_{V}\mid _{x}:=\nabla _{V}\mid _{x}+H^{o}_{V}\mid _{x},\qquad
(3.20^\prime )
\]
 where
\[
 H^{o}_{V}\mid _{x}
:=\bigl[\nabla _{V}|_{\gamma _{\mu (x)}}{\bigl(}
H^{\hbox{i.}}_{.j}(t,s^{o}_{\mu (x)};\gamma _{\mu (x)})
E_{i}|_{\gamma _{\mu (x)}}\otimes
\]
\[
 \otimes E^{j}|_{\gamma _{\mu (x)}}{\bigr)}\bigr]|_{t=s^{o}_{\mu
(x)}}.\qquad (3.21^{\prime\prime})
\]
 The so constructed derivation depends, of
course, on the initial $L$-transport along paths as well as on the used in
its definition family of paths $\{\gamma _{\lambda }\}$.

{\bf Proposition} ${\bf 3}{\bf .}{\bf 4}{\bf .} A$ given $L$-transport along
$\gamma $ in the tensor algebra over $M$ satisfies (3.8) and (3.17) iff the
generated by it map ${\cal D}^{\gamma }$, whose action is given by (2.10) or
(2.16), is a derivation of the tensor algebra over $\gamma (J)$, i.e. when
${\cal D}^{\gamma }_{s}$ is linear and
\[
 {\cal D}^{\gamma }_{s}(A\otimes B)=({\cal D}^{\gamma }_{s}A)\otimes
B+A\otimes ({\cal D}^{\gamma }_{s}B),\qquad (3.23)
\]
\[
{\cal D}^{\gamma }_{s}\circ C-C\circ {\cal D}^{\gamma }_{s}=0,\qquad (3.24)
\]
for arbitrary $C^{1}$tensor fields A and $B$ over $\gamma (J)$ and
contraction operator C.

{\bf Proof.} The linearity of ${\cal D}^{\gamma }_{s}$is a consequence from (2.11) or from the definitions of a derivation of the tensor algebra [3] and does not at all depend on the validity of (3.8) and (3.17).

If (3.8) and (3.17) are true, then with the help of (2.10) we see that from them follow (3.23) and (3.24) respectively. The same result is a consequence from the fact that in this case the considered $L$-transport in the tensor algebra over $M$ is an $S$-transport (see [4], definitions 2.1) for which, by proposition 3.5 from [4], the operator ${\cal D}^{\gamma }$is a derivation over $\gamma (J)$.

On the other hand, let ${\cal D}^{\gamma }$be a derivation of the mentioned tensor algebra. Then, by [4], proposition 3.5, there exists a unique $S$-transport generating ${\cal D}^{\gamma }$through (2.10). Therefore due to (2.14), the coefficients of this $S$-transport and of the considered $L$-transport coincide, which by [1], proposition 4.7 means these two transports along paths to coincide. But by definition (see [1], definition 2.1) any $S$-transport satisfies (3.8) and (3.17), consequently the investigated $L$-transport also satisfies them.\blacksquare

{\bf Proposition 3.5.} Every $S$-transport is an $L$-transport along (smooth) paths in the tensor algebra over M.

{\bf Proof.} This proposition follows from the comparison of definitions of $S$-transports $(cf. [4]$, definition $2.1), L$-transports $(cf. [1]$, definition 2.1) along paths, and $L$-transports along paths in the tensor algebra over a manifold (see (3.2) and (3.3)).\blacksquare

In the general case the "inverse proposition" to proposition 3.5 is not true. Because of this we shall consider the question when one $L$-transport along paths in the tensor algebra over $M$ is an $S$-transport along paths.

From the proof of proposition 3.4 one immediately derives the following two corollaries.

{\bf Corollary 3.4.} If the generated by (2.10) from a given $L$-transport along paths in the tensor algebra over $M$ operator ${\cal D}^{\gamma }$is a derivation of the tensor algebra over $\gamma (J)$, then this $L$-transport along $\gamma $ coincides with the $S$-transport along $\gamma $ which is generated from the defined by (3.20) derivation along $\gamma $.

{\bf Corollary} ${\bf 3}{\bf .}{\bf 5}{\bf .} A$ given $L$-transport along $\gamma $ in the tensor algebra over $M$ is an $S$-transport along $\gamma $ if and only if it satisfies the conditions (3.8) and (3.17) or, equivalently, iff the generated by it operator ${\cal D}^{\gamma }$is a derivation of the tensor algebra over $\gamma (J)$.

{\bf Corollary} ${\bf 3}{\bf .}{\bf 6}{\bf .} A$ given $L$-transport along paths in the tensor algebra over $M$ is (globally) an $S$-transport along paths iff the generated by it through (2.10) operator ${\cal D}^{\gamma }$is in fact a derivation of the tensor algebra along any path $\gamma $.

{\bf Proof.} This statement is a direct consequence from corollary 3.5 and (2.10).\blacksquare

\medskip
\medskip
 {\bf 4. SPECIAL BASES FOR LINEAR TRANSPORTS\\ ALONG PATHS IN TENSOR BUNDLES}

\medskip
For linear transports along paths in tensor bundles, of course, is valid proposition 3.1 of [1], according to which along any path there is a class of bases in which the transport's matrix is unit. But in the tensor bundles $(T^{p}_{.q}(M),\pi ,M)$ with $p+q\ge 1$ there exists a privilege set of bases, the one holonomic bases associated with different local coordinates. In this connection arises the question, which is a subject of the present section, when the described in the mentioned proposition bases are holonomic.

The next result shows that the $L$-transports along paths in the tangent and cotangent bundles over a manifold are Euclidean not only in a sense that they are such along any fixed path $\gamma  ($see [1], definition 3.1 and proposition 3.2), but also in a sense that along every part of $\gamma $ without selfintersections they are generated, through the described in [1], definition 2.4 way, from local holonomic bases, i.e. from local coordinates.

{\bf Proposition 4.1.} If $p+q=1$, then in $(T^{p}_{.q}(M),\pi ,M)$ for any $L$-transport $L^{\gamma }$along a path $\gamma :J  \to M$ without selfintersections there exists local coordinates in a neighborhood of (a part of$) \gamma (J)$ such that the matrix of $L^{\gamma }$is unit in the (field of) holonomic bases generated by them in $(T^{p}_{.q}(M),\pi ,M), p+q=1$ when they are restricted on the same neighborhood of $\gamma (J)$.

{\bf Proof.} By proposition 3.1 of [1] in the tangent (resp. cotangent) bundle of $M$ there exists a set of described in it (field of) bases, defined only on $\gamma (J)$, in which the matrix of $L^{\gamma }$is unit. By lemma 7 of [5] in a neighborhood of any part of $\gamma (J)$ lying in some coordinate neighborhood for any such basis there exist local coordinates for which the restriction on $\gamma (J)$ of the generated by them holonomic bases in the tangent (resp. cotangent) bundle coincide with the previous (field of) bases.

Consequently, in $(T^{p}_{.q}(M),\pi ,M), p+q=1$ any special for $L^{\gamma }$bases can be extended in a holonomic way on a neighborhood of
(a part of$) \gamma (J).\blacksquare $

If the path $\gamma $ has selfintersections, then, generally, along $\gamma $ there does not exist local coordinates with the described in [5], lemma 7 properties. The cause for this is that at the points of selfintersections, as a rule, the bases, in which the matrix of an $L$-transport is unit, are not uniquely defined or are not continuous. Therefore along any "piece" without selfintersection of an arbitrary path there exist local coordinates with the described properties and which are explicitly constructed in the proof of lemma 7 of [5]. But these coordinates admit continuation not far then the points of selfintersections, if any.

From the proof of proposition 4.1 also follows that in $(T^{p}_{.q}(M),\pi ,M), p+q=1$ any special for $L^{\gamma }$basis can be extended in a holonomic way outside (a part of$) \gamma (J)$ if $\gamma $ is without selfintersections. Evidently, nevertheless of the properties of $\gamma $ such an extension can be done also (and globally, i.e. on the whole set $\gamma (J))$ in an anholonomic way.

For the tensor bundles $(T^{p}_{.q}(M),\pi ,M)$ with $p+q\ge 2$ proposition 4.1 is generally not true. The only general exception of this are the $S$-transports, i.e. the $L$-transports along paths in the tensor algebra over $M$ consistent with the tensor product and commuting with the contractions. In fact, the matrix elements of these transports are given by (3.10) and $(3.18) (cf$. proposition 3.2). Therefore these matrices are unit matrices of the corresponding size in any special for the transport bases in the tangent (or cotangent) bundle over M. But the last bases can be chosen as holonomic ones (in a neighborhood (of a part) of any path; cf. proposition 4.1). Hence in the holonomic bases generated in this way in any
tensor space over $M$ the matrices of the $S$-transport are unit (of the corresponding size).

\medskip
\medskip
 {\bf 5. SECTIONS LINEARLY TRANSPORTED ALONG PATHS}

\medskip
Let a linear transport along paths $L$ in the vector bundle $(E,\pi ,B)$ be
given.

{\bf Definition 5.1.} The section $\sigma \in $Sec$(E,\pi ,B)$ is linearly
transported $(L$-transported), or undergoes an $L$-transport, along the path
$\gamma :J  \to B$ if for $t\in J$ and some $s\in J$, we have
\[
 \sigma (\gamma (t))=L^{\gamma }_{s  \to t}\sigma (\gamma (s)).\qquad (5.1)
\]
We say that $\sigma $ is $L$-transported if (5.1) holds for every $\gamma $.

{\bf Proposition 5.1.} If (5.1) holds for some $s\in J$, then it is true for every $s\in $J.

 {\bf Proof.} The result is a trivial corollary of (2.3).\blacksquare

{\bf Proposition 5.2.} The values of an $L$-transported (resp. along $\gamma
)$ section $\sigma $ are uniquely defined by fixing its value $\sigma
(x_{0})$ at an arbitrary given point $x_{0}\in B ($resp. $x_{0}\in \gamma
(J))$.

{\bf Proof.} The result follows from (5.1) for such $s$ for which $\gamma (s)=x_{0}.\blacksquare $

{ \bf Proposition 5.3.} The $C^{1}$section $\sigma $ is $L$-transported along
the path $\gamma $ iff it satisfies the equation
\[
 {\cal D}^{\gamma }\sigma =0,\qquad (5.2)
\]
where ${\cal D}^{\gamma }$ is defined from $L$ through (2.10).

{\bf Proof.} If $\sigma $ is $L$-transported along $\gamma $, then (5.2)
follows from (5.1) and (2.13).

On the opposite, let (5.2) holds. Fixing a basis $\{e_{i}(s)\}$
in $\pi^{-1}(\gamma (s)), s\in J$, we have $\sigma =\sigma ^{i}e_{i}$and
defining $\bar{\sigma }(s):=(\sigma ^{1}(\gamma (s)),\ldots  ,\sigma
^{\dim(\pi ^{-1_{(\gamma (s)))}}}(\gamma (s)))$, we see that, due to (2.16),
the  $eq. (5.2)$ is equivalent to
\[
\frac{d\bar{\sigma}(s)}{ds}+ \Gamma _{\gamma }(s)\bar{\sigma }(s)=0.
\qquad (5.3)
\]
 Substituting here $\Gamma _{\gamma }(s)=F^{-1}(s;\gamma
)dF(s;\gamma )/ds ($see (2.15)) and (2.9)), we get $d[F(s;\gamma )\bar{\sigma
}(s)]/ds=0$, i.e. $F(s;\gamma )\bar{\sigma }(s)=$const$=F(s_{0};\gamma
)\bar{\sigma }(s_{0})$ for a fixed $s_{0}\in J$ and, consequently
$\bar{\sigma }(s)=F^{-1}(s;\gamma )F(s_{0};\gamma )\bar{\sigma
}(s_{0})=H(s,s_{0};\gamma )\bar{\sigma }(s_{0})$ which, due to (2.2) and
(2.6), is equivalent to (5.1). So, $\sigma $ is $L$-transported along $\gamma
$ section.\blacksquare

{ \bf Proposition 5.4.} The maps (2.1) define an $L$-transport along $\gamma
$ if and only if for every $\sigma _{0}\in \pi ^{-1}(\gamma (s))$ the section
$\sigma \in $Sec$(\pi ^{-1}(\gamma (J),\pi \mid _{\gamma (J)},\gamma (J))$
defined by
\[
 \sigma (\gamma (t))=L^{\gamma }_{s  \to t}\sigma _{0}\qquad (5.4)
\]
 is a solution of the initial-value problem
\[
 {\cal D}^{\gamma }\sigma =0,
\quad \sigma (\gamma (s))=\sigma _{0},\qquad (5.5)
\]
where ${\cal D}^{\gamma }$is some derivation along $\gamma $, i.e. for it
(2.11) and (2.12) hold.

{\bf Proof.} If (2.1) defines an $L$-transport along $\gamma $, then by proposition 5.3, we have ${\cal D}^{\gamma }\sigma =0, {\cal D}^{\gamma }$being the defined from (2.10) derivation along $\gamma $, and $\sigma (\gamma (s))=\sigma _{0}$because of (2.4). On the contrary, let (5.5) holds. By proposition 4.7 of [1] there exists a unique $L$-transport $^\prime L^{\gamma }$along $\gamma $ generating ${\cal D}^{\gamma }$through (2.10) and having the same coefficients as ${\cal D}^{\gamma }$. Therefore, due to proposition 5.3 and definition 5.1, the unique solution of (5.5) is $\sigma (\gamma (t))=^\prime L^{\gamma }_{s  \to t}\sigma _{0}$. So (see(5.4)), we find $^\prime L^{\gamma }_{s  \to t}\sigma _{0}=L^{\gamma }_{s  \to t}\sigma _{0}$for every $\sigma _{0}$, i.e. $L^{\gamma }_{s  \to t}=^\prime L^{\gamma }_{s  \to t}$and hence the looked for $L$-transport along $\gamma $ coincides with
$  ^\prime L^{\gamma }.\blacksquare $

Let us note the evident corollary from the last proof that the coefficients of the $L$-transport along paths, entering in (5.4), and of ${\cal D}^{\gamma }$, appearing in (5.5), coincides and that these two operators generate each other in the way considered in [1].

 The last proposition gives us a ground to give

{\bf Definition 5.2.} The equality ${\cal D}^{\gamma }\sigma =0$ will be called an equation of the linear transport (the $L$-transport equation) along $\gamma $.

In any basis the matrix form of the equation of $L$-transports along $\gamma$
has the form (5.3) in which
$\Gamma_{\gamma }(s):= [\Gamma^{i}_{.j}(s;\gamma )]$ is the matrix of the
coefficients of some $L$-transport or a derivation along $\gamma $.

According to proposition 5.3 (or 5.4) any $L$-transported (resp. along $\gamma )$ section satisfies the equation of $L$-transports (resp. along $\gamma )$.

The special bases for an $L$-transport along $\gamma $ are characterized by
(see $(2.9), (2.15)$ and [1])
 \[
 H(t,s;\gamma )={\Bbb I}\ or\ \Gamma _{\gamma }(s)=0.\qquad (5.6)
\]
 So, in them (2.16) reduces to
\[
 {\cal D}^{\gamma }_{s}\sigma
=  \frac{d\sigma^i(\gamma(s))}{ds} e_{i}(s)\qquad (5.7)
\]
 and the equation of the $L$-transport
takes the trivial form $d\sigma ^{i}(\gamma (s))/ds=0$. Hence $\sigma $ is
$L$-transported along $\gamma $ iff in these bases $\sigma ^{i}=$const, a
fact which follows also directly from $(5.6), (2.6)$ and $(2.2): \sigma
^{i}(\gamma (t))=H^{i}_{.j}(t,s;\gamma )\sigma ^{j}(\gamma (s))=\sigma
^{i}(\gamma (s))=$const. Consequently in any special for an $L$-transport
basis along $\gamma $ the components of an $L$-transported sections are
constant along the path of transportation. In this sense the $L$-transported
sections of a vector bundles are analogous to the parallelly transported (or
(covariantly) constant) vector fields in an Euclidean space with respect to
Cartesian coordinates.

{\bf Proposition 5.5.} If the $L$-transport $L^{\gamma }$along $\gamma $ in tensor algebra over $M$ satisfies (3.8), then the function $f:\gamma (J)\to {\Bbb R}$ is $L$-transported along $\gamma $ iff it is a constant on $\gamma (J)$.

{\bf Proof.} If $f$ is $L$-transported along $\gamma  ($see definition 5.1), then $f(\gamma (t))=L^{\gamma }_{s\to t}f(\gamma (s))=f(\gamma (s))L^{\gamma }_{s\to t}(1)=f(\gamma (s))$ for any $s,t\in J$, i.e. $f(\gamma (s))=$const$=f(\gamma (s_{0}))$ for fixed $s_{0}\in J$ and every $s\in $J. On the opposite, if $f(\gamma (s))=c=$const$\in {\Bbb R}$, then $f(\gamma (t))=c=L^{\gamma }_{s\to t}(c)=L^{\gamma }_{s\to t}(f(\gamma (s)))$, i.e. $f$ is $L$-transported along $\gamma .\blacksquare $

\newpage

 {\bf 6. PATHS WITH LINEARLY TRANSPORTED\\ TANGENT VECTOR (L-PATHS)}

\medskip
Let a linear transport along paths $L$ in the tangent bundle $(T(M),\pi ,M)$ over a differentiable manifold $M$ be given. Below we sketch a scheme for an introduction of a class of paths in $M$ which, with respect to $L$, behave in the same way as the geodesics does with respect to the defining them parallel transport or a linear connection [3].

{\bf Definition 6.1.} A $C^{1}$path $\gamma :J  \to M$ is a path with a
linearly transported tangent vector, or simply an $L$-path, if its tangent
vector field  $\dot\gamma\in Sec(T(\gamma (J)),\pi \mid _{\gamma (J)},\gamma
(J))\subset  \subset (T(M),\pi .M)$ is $L$-transported along $\gamma $.

 By definition 5.1 the path $\gamma $ is an $L$-path iff
\[
\dot\gamma(t)=L^{\gamma }_{s\to t}\dot\gamma(s),\quad s,t\in J,\qquad (6.1)
\]
which, due to proposition 5.3, equivalently means that    satisfies the
$L$-transport equation along $\gamma $, i.e.
\[
 {\cal D}^{\gamma } \dot\gamma=0,\qquad (6.1^\prime )
\]
 where ${\cal D}^{\gamma }$ is given by (2.10).

If $L^{\gamma }_{s  \to t}$is a smooth transport, i.e. if it has $a C^{1}$dependence on $t$, then through any point $x\in M$ in any direction $X\in T_{x}(M)$ there is one and only one $L$-path. More precisely, it is true the following theorem which is an evident generalization of the corresponding theorem concerning geodesic paths in manifolds with affine connection (cf. e.g. [3]).

{\bf Theorem 6.1.} If $x\in M, X\in T_{x}(M), J$ is an ${\Bbb R}$-interval
and $s_{0}\in J$ is fixed, then there exist a unique $L$-path $\gamma :J  \to
M$, such that
\[
\gamma (s_{0})=x,\quad   \dot\gamma(s_{0})=X.   \qquad (6.2)
\]

 {\bf Proof.} From (6.1) for $s=s_{0}$ and (6.2), we see that the statement
of the theorem is equivalent to the existence of a unique path $\gamma $
having the properties
\[
\dot\gamma(t)=L^{\gamma }_{s_{0}\to t}X,\qquad (6.3a)
\]
\[
 \gamma (s_{0})=x . \qquad (6.3b)
\]
Due to (2.2) and (2.6) in local coordinates (6.3a) reduces to a first
order system of ordinary differential equations with respect to the local
coordinates of $\gamma (t)$ which, due to the initial condition (6.3b), in
accordance with the conditions of the theorem  and the theorems for existence
and uniqueness of such systems [6] has a unique solution $\gamma :J  \to
$M.\blacksquare

Let us write the initial-value problem (6.3) in an equivalent but more convenient from practical view-point form, which is near to that in a case of geodesic paths [3].

Let ${\cal D}^{\gamma }$be the generated from the given $L$-transport $L$
along $\gamma $ derivation (see (2.10)). Due to proposition 5.4 the
initial-value problem (6.3) is equivalent to
\[
 {\cal D}^{\gamma }(\dot\gamma)=0,\qquad (6.4a)
\]
\[
\dot\gamma(s_{0})=X, \gamma(s_{0})=x,\qquad (6.4b)
\]
i.e. satisfies the $L$-transport equation along  $\gamma $ under the initial
conditions (6.2).

If in some local basis the transport $L$ is given by its coefficients
$\Gamma ^{i}_{.j}(s;\gamma ) ($see (2.15)), then in it, according to (5.3),
the equation (6.4a) takes a form analogous to that of the canonical geodesic
equation [3]:
\[
\frac{d\dot\gamma^i(s)}{ds} + \Gamma ^{i}_{.j}(s;\gamma ) \dot\gamma^{j}(s)=0,
\quad s\in J.  \qquad (6.5)
\]

As a consequence of theorem 6.1 the equation  (6.4a) or the system (6.5) can
be called {\it equation} or {\it a system of equations of the} $L${\it
-paths.}

Evidently $(cf. [3])$, the $L$-paths generalize the concept of geodesic paths (curves) to which they reduce when the transport $L$ is a parallel transport corresponding to a covariant differentiation (linear connection$) \nabla $ or, equivalently, when ${\cal D}^{\gamma }$is a covariant differentiation along $\gamma $, i.$e {\cal D}^{\gamma }=\nabla _{\cdot }$for a covariant differentiation $\nabla $.

{\bf Proposition 6.1.} Along any $L$-path there exist (a class of) local holonomic bases in which it is defined as a linear function of its parameter.

{\bf Proof.} Let us consider any special for $L$ basis along $\gamma $. In
it (5.6) holds, i.$e (6.5)$ reduces to
\[
\frac{d\dot\gamma^i(s)}{ds} =0.\qquad (6.6)
\]
By lemma 7 from [5] locally, i.e. in a neighborhood of any part of $\gamma
(J)$ lying in only one coordinate neighborhood in which $\gamma $ is without
selfintersections, this basis can be extended in a holonomic way outside
$\gamma (J)$. So, there are local coordinates $\{x^{i}\}$ in which
$\dot\gamma^{i}(s)=d\gamma ^{i}(s)/ds$ and also (6.6) are true.
Therefore, we have
\[
\frac{d^2\gamma^i(s)}{ds^2} =0,\qquad (6.7)
\]
 the general solution of which is
\[
 \gamma ^{i}(s)=X^{i}(s-s_{0})+x^{i}\qquad (6.8)
\]
 for some constants $s_{0}\in J, X^{i}$and $x^{i}.\blacksquare $

Comparing (6.8) (see the last proof) and (6.4b) we see $x^{i}$and $X^{i}$to be, respectively, the coordinates of the point $\gamma (s_{0})$ and the components of the vector   $(s_{0})$ at it in the considered special holonomic basis.

\medskip
\medskip
 {\bf 7. CONCLUSION}

\medskip
Here we have considered only a few examples of usage of linear transports along paths in vector bundles. Some of them, in particular the theory of $L$-paths, as well as the applications of the $L$-transport along paths to physical problems will be investigated in details elsewhere.

At the end we want to  make a comment on the special bases for a linear transport along paths in the tangent
bundles over a manifold when it is a parallel transport associated to a linear connection with local coefficients $\Gamma ^{i}_{.jk}$.

In this case the transport's coefficients along $\gamma $ are (see [1],
Sect. 5)
\[
 \Gamma ^{i}_{.j}(s;\gamma )
=\Gamma ^{i}_{.jk}(\gamma (s)) \dot\gamma^{k}(s).\qquad (7.1)
\]

If $\{E_{i^\prime }\}$ is a special along $\gamma $ for the transport basis
(see [1], Sect. 3), then $\Gamma ^{i^\prime }_{..j^\prime }(s;\gamma )=0$,
i.e.
\[
\Gamma ^{i^\prime }_{..j^\prime k^\prime }(\gamma (s))
\dot\gamma^{k^\prime}(s)=0.\qquad (7.2)
\]

 As $\{E_{i^\prime }\}$ itself depends, generally, on
$\gamma $, from here one can not conclude that $\Gamma ^{i^\prime
}_{..j^\prime k^\prime }(\gamma (s))=0$. But in [5], corollary 11 we proved
the existence of a class of local bases, defined in a neighborhood of $\gamma
(J)$, in any one of which the connection's components vanish on $\gamma (J)$.
Evidently (see (7.1)), these bases are special for the corresponding to the
connection parallel transport. Comparing the arbitrariness in the definitions
of the bases belonging to the considered two sets of special bases, for the
connection (see corollary 11 and proposition 2 from [5]) and for assigned to
it parallel transport (see proposition 3.1 from [1]), we conclude these two
sets to be identical (on $\gamma (J))$.

 Hence, for a linear connection on the set $\gamma (J)$, defined by a path
$\gamma :J  \to M$, in any basis in which the connection's coefficients
vanish also vanish the coefficients of the corresponding to it parallel
transport and vice versa.

\medskip
\medskip
 {\bf ACKNOWLEDGEMENTS}

\medskip
 The author expresses his gratitude to Prof. Vl. Aleksandrov (Institute of
Mathematics of Bulgarian  Academy of Sciences) for constant interest in this
work and stimulating discussions.

This research was partially supported by the Fund for Scientific Research of
Bulgaria under contract Grant No. $F 103$.

\medskip
\medskip
 {\bf REFERENCES}

\medskip
1.  Iliev  B.Z.,  Linear  transports  along  paths  in  vector bundles. I. General theory, Communication JINR, $E5-93-239$, Dubna, 1993. \par
2.  Iliev B.Z., The equivalence principle in spaces with a linear transport along paths, Proceedings of the 5-th seminar "Gravitational energy and gravitational waves", JINR, Dubna, $16-18$ may 1992, Dubna, $1993, pp.69-84 ($In Russian).\par
3.  Kobayashi S., K. Nomizu, Foundations of Differential Geometry, Vol. 1, Interscience Publishers, New York-London, 1963.\par
4.  Iliev B.Z., Parallel transports in tensor spaces generated by derivations of tensor algebras, Communication JINR, $E5-93-1$, Dubna, 1993.\par
5.  Iliev B.Z., Special bases for derivations of tensor algebras. II. Case along paths, Communication JINR, $E5-92-508$, Dubna, 1992.
\par
6.  Hartman Ph., Ordinary Differential Equations, John Wiley \& Sons, New York-London-Sydney, $1964, ch$.IV, \S1.

\newpage
\vspace*{6ex}
\noindent
 Iliev B. Z.\\[5ex]

\noindent
 Linear Transports along Paths in Vector Bundles\\
 II. Some Applications\\[5ex]

\medskip
\medskip
 The linear transports along paths in vector bundles introduced in Ref. [1]
are applied to the special case of tensor bundles over a given differentiable
manifold. The ties with the transports along paths generated by derivations
of tensor algebras are investigated. A possible generalization of the theory
of geodesics is proposed when the parallel transport generated by a linear
connection is replaced with an arbitrary linear transport along paths in the
tangent bundle.\\[5ex]

\medskip
\medskip
The investigation has been performed at the Laboratory of Theoretical
Physics, JINR.

\end{document}